%% file: post-simmer-quasi-unitary-equivalence-arxiv.tex
\documentclass[12pt,reqno,a4paper]{amsart}            

\usepackage{hyperref}

\usepackage{scrlayer-scrpage}
\usepackage{scrtime}     

\usepackage{amsmath, amssymb}
\usepackage{amsthm}

\usepackage{accents}     
\usepackage{mathtools}  

\usepackage{bbm}         

\usepackage{mathrsfs}    

\usepackage{graphicx}                

\usepackage[totalwidth=450pt, totalheight=750pt]{geometry}
\input{\localpath olaf-ams-layout-macros.tex}
\input{\localpath olaf-thm-environment-macros.tex}
\input{\localpath olaf-refs-unified.tex}

\input{\localpath olaf-math-symbols.tex}

\input{\localpath olaf-other-macros-myenumi-myparagraph.tex}

\input{\localpath olaf-spaces-macros.tex}

\providecommand{\myfont}{}
\usepackage{hyperref}
\usepackage{nicefrac}
\usepackage{enumerate}

\newcommand{\deltaA}{\delta_{\mathrm a}}

\newcommand{\deltaC}{\delta_{\mathrm c}}

\newcommand{\Cellreg}{C_{\mathrm{ell.reg}}}
\newcommand{\Cext}{C_{\mathrm{ext}}}

\begin{document}

\title[Quasi-unitary equivalence and generalised norm resolvent
convergence] {Quasi-unitary equivalence and generalised norm resolvent
  convergence}

\author{Olaf Post}
\address{Fachbereich 4 -- Mathematik,
  Universit\"at Trier,
  54286 Trier, Germany}
\email{olaf.post@uni-trier.de}

\author{Jan Simmer}%
\address{Fachbereich 4 -- Mathematik,
  Universit\"at Trier,
  54286 Trier, Germany (till 2020)}
\email{jan.simmer@yahoo.de}

\ifthenelse{\isundefined \finalVersion}
           {\date{\today, \thistime,  \emph{File:} \texttt{\jobname.tex}}}
           {\date{\today}}  

\begin{abstract}
  The purpose of this article is to give a short introduction to the
  concept of quasi-unitary equivalence of quadratic forms and its
  consequences.  In particular, we improve an estimate concerning the
  transitivity of quasi-unitary equivalence for forms.  We illustrate
  the abstract setting by two classes of examples.
\end{abstract}


\maketitle



%
%


%
\section{Generalised norm-resolvent convergence}
\label{sec:que}
%

In this article, we give an overview of the concept of
\emph{quasi-unitary equivalence} for non-negative and self-adjoint
operators and quadratic forms acting in \emph{different} Hilbert
spaces. Quasi-unitary equivalence provides a sort of ``distance''
between two operators resp.\ two quadratic forms.  The concept was
introduced by the first author in~\cite{post:06} and later explained
in great detail in the monograph~\cite[Ch.~4]{post:12}.  In this
article, we also improve an estimate concerning the transitivity of
the notion of quasi-unitary equivalence of quadratic forms in
\Prp{trans.q-u-e} (cf.~\cite[Prp.~4.4.16]{post:12}).  Moreover, we
also specify the convergence rate of functions of the operators
in~\Thm{func.calc} which also follows from quasi-unitary equivalence
and thus simplifying some earlier results (cf.~\cite[Thm.~4.2.9 and
following pages]{post:12}).  In particular, we make both estimates
more explicit.

To illustrate the strength of the concept, we give several examples.
We start with the approximation of Laplacians on post-critically
finite fractals such as the Sierpi\'nski gasket by their
finite-dimensional analogues, see \Sec{pcf-fractals}; more details can
be found in~\cite{post-simmer:18}.
In~\cite{post-simmer:20,post-simmer:pre18c} we extend the results
to \emph{magnetic} Laplacians and more general spaces such as finitely
ramified fractals.  In~\cite{post-simmer:21} we use a similar
strategy to compare Laplacians on discrete graphs with metric spaces
such as metric graphs and graph-like manifolds.  In particular, using
the above-mentioned transitivity, we could show that Laplacians on a
post-critically finite fractal can be approximated by a sequence of
Laplacians on manifolds.

The concept of quasi-unitary equivalence is not restricted to such
discretisations: it can also be applied in other situations where the
underlying spaces change: it was originally developed for a family of
thin manifolds shrinking to a so-called metric (or quantum) graph
(see~\cite{post:06}), but can also be applied to other (drastic)
changes of the space as in~\cite{anne-post:21}.  We give a flavour
of such arguments in \Sec{neu-obstacles} where we apply the concept to
the case of a manifold with (small) obstacles taken out: we show that
the Neumann Laplacian on the remaining set is close to the original
(unperturbed) Laplacian.

Other applications are possible, e.g.~in~\cite{khrabustovskyi-post:18}
we could apply it to homogenisation problems and show a generalised
\emph{norm} resolvent convergence.  Typically, results in
homogenisation theory only include \emph{strong} resolvent
convergence.  Due to the lack of spectral exactness (the limit
spectrum can suddenly shrink: there may be approximating sequences
that do not correspond to spectral values in the limit, so-called
\emph{spectral pollution}).  Norm resolvent convergence (and also our
generalised version of it, see \Def{gnrs}) implies that the spectra
converge, see \Subsec{que-consequences} for details.  Our concept is
also closely related to generalised norm resolvent convergence in the
sense of Weidmann, see~\cite[Sec.~9.3]{weidmann:00}.  We also refer to
the recent work~\cite{boegli:18} and references therein; we deal with
these aspects in a forthcoming publication.

\subsection{Quasi-unitary equivalence for operators}
\label{ssec:operator-que}

We first start defining a ``distance'' between two non-negative and
self-adjoint operators $\Delta$ and $\wt\Delta$ acting in
\emph{different} Hilbert spaces $\HS$ and $\wt\HS$.  The distance is
expressed in terms of a parameter $\delta \ge 0$, and appears in the
concept of \emph{$\delta$-quasi-unitary equivalence}, which we will
explain now.

Associated with such a $\Delta$ we define a so-called \emph{scale of
  Hilbert spaces} $\HS^k:=\dom \Delta^{k/2}$ with norm
$\norm[k] f :=\norm[\HS]{(\Delta+1)^{k/2}f}$ for $k \ge 0$.  For
negative powers, we let $\HS^{-k}$ be the completion of $\HS$ under
the norm $\norm[-k] f := \norm[\HS]{(\Delta+1)^{-k/2}f}$; moreover the
inner product $\iprod[\HS]\cdot \cdot$ extends continuously onto the
dual pairing $\HS^{-k} \times \HS^k$.  Similarly, we have a scale of
Hilbert spaces $\wt \HS^k$ associated with $\wt \Delta$.
\begin{definition}
  \label{def:quasi-uni.op}
  Let $\delta \ge 0$.
  \begin{subequations}
    \begin{enumerate}
    \item We say that a linear operator $\map J \HS {\wt\HS}$ is
      \emph{$\delta$-quasi-unitary} with \emph{$\delta$-quasi-adjoint}
      $\map {J'} {\wt \HS}\HS$ (for the operators $\Delta$ and $\wt
      \Delta$) if
      \begin{gather}
        \label{eq:quasi-uni.op.a}
        \norm{Jf}\le (1+\delta) \norm f, \quad
        \bigabs{\iprod {J f} u - \iprod f {J' u}}
        \le \delta \norm f \norm u
        \qquad (f \in \HS, u \in \wt \HS),\\
        \label{eq:quasi-uni.op.b}
        \norm{f - J'Jf}
        \le \delta \norm[2] f, \quad
        \norm{u - JJ'u}
        \le \delta \norm[2] u \qquad (f \in \HS^2, u \in \wt \HS^2).
      \end{gather}
      We call $J$ and $J'$ \emph{identification operators}.

    \item We say that the operators $\Delta$ and $\wt \Delta$ are
      \emph{$\delta$-close} if
      \addtocounter{equation}{1}
      \begin{equation}
        \label{eq:quasi-uni.op.d}
        \bigabs{\iprod[\wt \HS]{Jf}{\wt \Delta u} -
                \iprod[\HS] {J\Delta f}u}
        \le \delta \norm[2] f \norm[2] u
        \qquad (f \in \HS^2, u \in \wt \HS^2).
      \end{equation}

    \item We say that $\Delta$ and $\wt \Delta$ are
      \emph{$\delta$-(operator-)quasi-unitarily equivalent},
      if~\eqref{eq:quasi-uni.op.a}--\eqref{eq:quasi-uni.op.d} are
      fulfilled, i.e., we have the following equivalent operator norm
      estimates
      \begin{gather}
        \label{eq:quasi-uni.op.a'}
        \tag{\ref{eq:quasi-uni.op.a}'}
        \norm J \le 1+\delta, \qquad \norm{J^* - J'} \le \delta\\
        \label{eq:quasi-uni.op.b'}
        \tag{\ref{eq:quasi-uni.op.b}'}
        \norm{(\id_\HS - J'J)R} \le \delta, \qquad
        \norm{(\id_{\wt \HS} - J J')\wt R} \le \delta,\\
        \label{eq:quasi-uni.op.d'}
        \tag{\ref{eq:quasi-uni.op.d}'}
        \norm{\wt R J - J R} \le \delta,
      \end{gather}
      where $R:=(\Delta+1)^{-1}$ and $\wt R:=(\wt \Delta + 1)^{-1}$.
    \end{enumerate}
  \end{subequations}
\end{definition}
Note that we also have
\begin{equation}
  \label{eq:quasi-uni.op.a''}
  \norm{J'}
  \le \norm{J'-J^*} + \norm{J^*}
  \le 1+2\delta,
\end{equation}
using $\norm {J^*}$ and~\eqref{eq:quasi-uni.op.a'}.

Note that if $\delta=0$ in the above definition then  $J$ is unitary
with inverse $J^*=J^\prime$ by~\eqref{eq:quasi-uni.op.a'}
and~\eqref{eq:quasi-uni.op.b'}. Thus the corresponding
operators $\Delta$ and $\wt\Delta$ are \emph{unitarily equivalent}
by~\eqref{eq:quasi-uni.op.d'}. Hence, quasi-unitary equivalence
\emph{generalises unitary equivalence}.

\begin{definition}
  \label{def:gnrs}
  Let $\Delta_m$ be a self-adjoint and non-negative operator acting in
  $\HS_m$ for $m \in \bar \N := \N \cup \{\infty\}$.  We say that the sequence
  $\{\Delta_m\}_{m\in\N}$ \emph{converges in generalised norm resolvent sense
    (with error estimate $\{\delta_m\}_{m\in\N}$)} to $\Delta_\infty$, if
  $\Delta_m$ and $\Delta_\infty$ are $\delta_m$-quasi-unitarily
  equivalent with $\delta_m \to 0$.
\end{definition}
The notion also \emph{generalises} the concept of \emph{norm resolvent
  convergence}: Assume that the operators all act in the same Hilbert
space, i.e., that $\HS:=\HS_m=\HS_\infty$ for all $m\in\N$.  Then the
sequence $\{\Delta_m\}_{m\in\N}$ converges in \emph{norm resolvent
  sense} to $\Delta_\infty$ if
\begin{equation}
  \label{eq:nrs}
  \bignorm{(\Delta_m+1)^{-1} - (\Delta_\infty +1)^{-1}} \to 0
  \qquad\text{as}\quad m\to\infty.
\end{equation}
If we choose $J$ and $J'$ to be the identity operator on $\HS$,
then~\eqref{eq:quasi-uni.op.a'} and~\eqref{eq:quasi-uni.op.b'} are
fulfilled with $\delta_m=0$, and~\eqref{eq:quasi-uni.op.d'} with
$\delta_m \to 0$ is equivalent with~\eqref{eq:nrs}.  A sequence of
operators also converges in \emph{generalised} norm resolvent sense if
there is a sequence of unitary operators
$\map {J_m}{\HS_m}{\HS_\infty}$ such that
$\norm{J_m(\Delta_m+1)^{-1}J^*_m-(\Delta_\infty+1)^{-1}} \to 0$.

The notion of operator-quasi-unitary equivalence is transitive in the
following sense (the proof is similarly to the one of
\Prp{trans.q-u-e}, and we slightly improved the error term
$\hat \delta$ compared to the one given in~\cite[Prp.~4.2.5]{post:12}):
\begin{proposition}[{\cite[Prp.~4.2.5]{post:12}}]
  \label{prp:trans.q-u-e.op}
  Assume that $\delta,\wt \delta \in [0,1]$.  Assume in addition that
  $\Delta$ and $\wt \Delta$ are $\delta$-quasi-unitarily equivalent
  with identification operators $J$ and $J'$, and that $\wt \Delta$
  and $\widehat \Delta$ are $\wt \delta$-quasi-unitarily equivalent
  with identification operators $\wt J$ and $\wt J'$.  Then $\Delta$
  and $\widehat \Delta$ are $\hat \delta$-quasi-unitarily equivalent
  with identification operators $\hat J=\wt J J$ and
  $\hat J'=J'\wt J'$, where $\hat \delta=5\delta+5\wt\delta$.
\end{proposition}

\subsection{Quasi-unitary equivalence for energy forms}
\label{ssec:que-energy}

It is actually more convenient to start with the quadratic forms
$\energy$ and $\wt \energy$ associated with the non-negative operators
$\Delta$ and $\wt \Delta$, and develop a slightly more elaborated
version of quasi-unitary equivalence.  This approach avoids dealing
with the sometimes complicated operator domains and graph norms.
Nevertheless, in applications, the more elaborated conditions are
easily verified.

Let $\HS$ and $\wt\HS$ be two separable (complex) Hilbert spaces. We
say that $\energy$ is an \emph{energy form in $\HS$} if $\energy$ is a
closed, non-negative quadratic form in $\HS$, i.e., if $\energy(f)
\coloneqq \energy(f,f)$ for some sesquilinear form $\map
{\energy}{\HS^1 \times\HS^1} \C$, denoted by the same symbol, if
$\energy(f)\ge 0$ and if $\HS^1:=\dom \energy$, endowed with the norm
defined by
\begin{equation}
  \label{eq:qf.norm}
  \normsqr[\energy] f
  \coloneqq \normsqr[\HS] f + \energy(f),
\end{equation}
is itself a Hilbert space and dense (as a set) in $\HS$.  Note that
$\norm[1] f = \norm[\energy] f$ and that
$\norm[\energy] f \le \norm[2]f$ in the terminology of
\Subsec{operator-que}.  We call the corresponding non-negative, self
adjoint operator $\Delta$ (see e.g.~\cite[Sec.~VI.2]{kato:66}) the
\emph{Laplacian} associated with $\energy$.  Similarly, let
$\wt \energy$ be an energy form in $\wt \HS$ with Laplacian
$\wt \Delta$.

We now also need identification operators $J^1$ and $J^{\prime 1}$
acting on the form domains.
\begin{definition}
  \label{def:quasi-uni}
  \begin{subequations}
    \label{eq:quasi-uni}
    Let $\delta \ge 0$, $\map J \HS {\wt \HS}$ and $\map {J'}{\wt
      \HS}\HS$, resp.\ $\map {J^1} {\HS^1} {\wt \HS^1}$ and
    $\map{J^{\prime1}} {\wt \HS^1}{\HS^1}$ be linear operators on the
    Hilbert spaces and energy form domains.
    \begin{enumerate}
    \item We say that $J$ is \emph{$\delta$-quasi-unitary} with
      \emph{$\delta$-quasi-adjoint} $J'$ (for the energy forms
      $\energy$ and $\wt \energy$) if
      \begin{align}
        \label{eq:quasi-uni.a}
        &\norm{Jf}\le (1+\delta) \norm f,
        &\bigabs{\iprod {J f} u - \iprod f {J' u}}
        &\le \delta \norm f \norm u
        &(f \in \HS, u \in \wt \HS),\\
        \label{eq:quasi-uni.b}
        &\norm{f - J'Jf}
        \le \delta \norm[\energy] f,
        &\norm{u - JJ'u}
        &\le \delta \norm[\wt \energy] u
        &(f \in \HS^1, u \in \wt \HS^1).
      \end{align}

    \item We say that $J^1$ and $J^{\prime1}$ are
      \emph{$\delta$-compatible} (with identification operators $J$
      and $J'$) if
      \begin{equation}
        \label{eq:quasi-uni.c}
        \norm{J^1f - Jf}\le \delta \norm[\energy]f, \quad
        \norm{J^{\prime1}u - J'u} \le \delta \norm[\wt \energy] u \qquad
        (f \in \HS^1, u \in \wt \HS^1).
      \end{equation}

    \item We say that the energy forms $\energy$ and $\wt \energy$ are
      \emph{$\delta$-close} if
      \begin{equation}
        \label{eq:quasi-uni.d}
        \bigabs{\wt \energy(J^1f, u) - \energy(f, J^{\prime1}u)}
        \le \delta \norm[\energy] f \norm[\wt \energy] u \qquad
        (f \in \HS^1, u \in \wt \HS^1).
      \end{equation}

    \item We say that $\energy$ and $\wt \energy$ are
      \emph{$\delta$-quasi-unitarily equivalent},
      if~\eqref{eq:quasi-uni.a}--\eqref{eq:quasi-uni.d} are fulfilled.
    \end{enumerate}
  \end{subequations}
\end{definition}


We have the following relation between quasi-unitary equivalence for
quadratic forms and operators; the last conclusion has already been
shown in~\cite[Prp.~4.4.15]{post:12}:
\begin{proposition}
  \label{prp:q-u-e.qf.op}
  If the forms $\energy$ and $\wt \energy$ are
  $\delta$-quasi-unitarily equivalent then we have
  \begin{equation}
    \label{eq:res.z}
    \norm{\wt R(z)J-JR(z)} \le C(z) \delta,
  \end{equation}
  where $R(z):=(\Delta-z)^{-1}$ and $\wt R(z):=(\wt \Delta-z)^{-1}$
  for $z \in \C \setminus (\spec \Delta \cup \spec {\wt \Delta})$ and
  \begin{equation*}
    C(z)
    := 4\Bigl(1 + \frac{\abs{z+1}}{d(z,\spec \Delta \cup \spec{\wt \Delta})}
    \Bigr)^2.
  \end{equation*}
  In particular, the associated operators $\Delta$ and $\wt \Delta$
  are $4\delta$-quasi-unitarily equivalent.
\end{proposition}
\begin{proof}
  For $g \in \HS$ and $v \in \wt \HS$, we have
  \begin{align*}
    &\bigabs{\iprod{(\wt R(z)J-JR(z))g} v}
    = \bigabs{\iprod g {J^*\wt R(\conj z)v}-\iprod{JR(z)g} v}\\
    \quad&
    = \bigabs{\iprod{\Delta f}{J^*u}-\iprod{Jf}{\wt \Delta u}}\\
    \quad&
    \le \bigabs{\bigiprod{\Delta f}{\bigl((J^*-J')+(J'-J^{\prime 1})\bigr) u}}
      +  \bigabs{\energy(f,J^{\prime 1}u)-\wt \energy(J^1 f,u)}
      +  \bigabs{\iprod{(J^1-J)f}{\wt \Delta u}}\\
    \quad &
    \le 2\delta \norm[2]f \norm u + \delta \norm[1]f\norm[1]u
       + \delta \norm f \norm[2] u
    \le 4\delta \norm[2]f \norm[2] u,
  \end{align*}
  where $f=R(z)g$ and $u=\wt R(\conj z) v$.  We have
  $\norm[2] f = \norm{(\Delta+1)R(z)g} \le \norm{(\Delta+1)R(z)} \norm
  g$ and
  \begin{equation*}
    \norm{(\Delta+1)R(z)}
    = \sup_{\lambda \in \spec \Delta} \frac{\lambda+1}{\abs{\lambda-z}}
    \le 1+ \sup_{\lambda \in \spec \Delta} \frac{\abs{z+1}}{\abs{\lambda-z}}
    = 1+ \frac{\abs{z+1}}{d(z,\spec \Delta)}
  \end{equation*}
  using the spectral theorem.  A similar estimate holds for
  $\norm[2] u$.  In particular, the resolvent estimate follows.  For
  the second statement, note that for $z=-1$ we have $C(-1)=1$,
  hence~\eqref{eq:quasi-uni.op.d} holds with $4\delta$.  The remaining
  estimates~\eqref{eq:quasi-uni.op.a}--\eqref{eq:quasi-uni.op.b}
  follow from the quasi-unitary equivalence of the forms and the fact
  that $\norm[1] f \le \norm[2] f$ and similarly for $u$.
\end{proof}
In particular, if we choose the rough estimate
$\spec \Delta \cup \spec {\wt \Delta} \subset [0,\infty)$, then
$C(z) \le 1 + \abs{z+1}/d(z,[0,\infty))$.  For $\Re z \ge 0$, the
latter equals $1+\abs{z+1}/\abs{\Im z}$ and for $\Re z<0$ the latter
equals $1+\abs{z+1}/\abs z$.  Hence, we have
\begin{equation*}
  C(z) \le 1 + \frac{\abs{z+1}}{\abs{\Im z}}
  \quadtext{resp.}
  C(z) \le 1 + \frac{\abs{z+1}}{\abs z}
\end{equation*}
for $\Re z \ge 0$ resp.\ $\Re z < 0$.

Let us mention a special case here, namely $\delta=0$ in
\eqref{eq:quasi-uni.a}--\eqref{eq:quasi-uni.c}.  In this situation,
$J$ is a unitary operator with $J'=J^*$, $J^1=J \restr {\HS^1}$ and
$J^{\prime 1}=J^* \restr {\wt \HS^1}$; hence without loss of
generality we can assume $\HS=\wt \HS$, $J=J'=\id_\HS$ and
$\dom \energy = \dom \wt \energy$.  In particular, $\energy$ and
$\wt \energy$ are $\delta$-quasi-unitarily-equivalent if and only if
\begin{subequations}
  \begin{equation}
    \label{eq:quasi-uni.d'}
    \abs{\wt\energy(f,u)-\energy(f,u)}
    \le \delta \norm[\energy] f \norm[\wt \energy] u
  \end{equation}
  for all $f,u \in \HS^1:=\dom \energy = \dom \wt \energy$. Using the
  fact that $\energy$ and $\wt \energy$ are symmetric, it is
  sufficient if~\eqref{eq:quasi-uni.d'} only holds for $f=u$, i.e.,
 ~\eqref{eq:quasi-uni.d'} is equivalent with
  \begin{equation}
    \label{eq:quasi-uni.d''}
    \abs{\wt\energy(f)-\energy(f)}
    \le \hat \delta \normsqr[\energy] f
  \end{equation}
\end{subequations}
for all $f \in \HS^1$.  For the
implication~\eqref{eq:quasi-uni.d'}$\Rightarrow$\eqref{eq:quasi-uni.d''}
one can use $\hat \delta=\delta \sqrt{(2+\delta)/(2-\delta)}$
(provided $\delta < 2$) and
for~\eqref{eq:quasi-uni.d''}$\Rightarrow$\eqref{eq:quasi-uni.d'} one
can use $\delta=\hat \delta/\sqrt{1-\hat \delta}$ (provided
$\hat \delta < 1$).  This situation has also been studied
in~\cite{brasche-fulsche:18}; basically, their Theorem~2 is the
implication
\eqref{eq:quasi-uni.d''}$\Rightarrow$\eqref{eq:quasi-uni.d'} together
with \Prp{q-u-e.qf.op} (with $z=-1$ and $4\delta$ replaced by
$\delta$, as $\delta=0$
in~\eqref{eq:quasi-uni.a}--\eqref{eq:quasi-uni.c}).

In particular, if $\{\energy_m\}_{m\in\N}$ is a sequence of energy
forms acting in the \emph{same} Hilbert space as $\energy_\infty$,
i.e., $\HS_m=\HS_\infty$ with the \emph{same} domain
$\dom\energy_m=\dom\energy_\infty$ for all $m\in\N$, then (with all
identification operators being the corresponding identity
operators)~\eqref{eq:quasi-uni.a}--\eqref{eq:quasi-uni.c} are trivially
fulfilled with $\delta=0$.  Moreover,~\eqref{eq:quasi-uni.d} with
$\delta=\delta_m \to 0$ is equivalent with
\begin{equation*}
  \bigabs{\energy_\infty(f)-\energy_m(f)}
  \le \hat \delta_m \normsqr[\energy_\infty] f
\end{equation*}
for all $f \in \dom \energy_\infty$ with $\delta_m$ and
$\hat \delta_m$ related as above.  This is the classical situation of
Kato~\cite[Thm.~VI.3.6]{kato:66}
or~\cite[Thm.~VIII.25(c)]{reed-simon-1}, and we conclude (using
\Prp{q-u-e.qf.op}) that the operators $\Delta_m$ associated with
$\energy_m$ \emph{converge} to $\Delta_\infty$ in \emph{norm resolvent
  sense}, see~\eqref{eq:nrs}.  Note that both classical results do not
state the convergence speed of the norm of the resolvent difference.

Another useful implication is the transitivity of quasi-unitary
equivalence for energy forms; it was originally proved
in~\cite[Prp.~4.4.16]{post:12}; we give here a simpler proof.
\begin{proposition}
  \label{prp:trans.q-u-e}
  \begin{subequations}
    Let $\delta, \wt \delta \in [0,1]$.  Assume that $\energy$ and
    $\wt \energy$ are $\delta$-quasi-unitarily equivalent with
    identification operators $J$, $J^1$, $J^\prime$ and
    $J^{\prime 1}$.  Moreover, assume that $\wt \energy$ and
    $\widehat \energy$ are $\wt\delta$-quasi-unitarily equivalent with
    identification operators $\wt J$, $\wt J^1$, $\wt J^\prime$ and
    $\wt J^{\prime 1}$. Assume in addition that, for all $f\in\HS^1$
    and $w\in\widehat\HS^1$,
    \begin{equation*}
      \norm[\wt\energy]{J^1f}
      \leq (1+\delta)\norm[\energy]{f}
      \quadtext{and}
      \norm[\wt\energy]{\wt J^{\prime 1}w}
      \leq (1+\wt \delta)\norm[\widehat\energy]{w}.
    \end{equation*}
    Then $\energy$ and $\widehat \energy$ are
    $\hat \delta$-quasi-unitarily equivalent with
    $\hat \delta=14(\delta+\wt \delta)$.
  \end{subequations}
\end{proposition}

\begin{proof}
  We define the identification operators by $\widehat J:=\wt JJ$,
  $\widehat J^1:=\wt J^1 J^1$,
  $\widehat J^\prime:=J^\prime\wt J^\prime$ and
  $\widehat J^{\prime 1}:=J^{\prime 1}\wt J^{\prime 1}$ and we set
  $R:=(\Delta+1)^{-1}$, $\wt R:=(\wt\Delta+1)^{-1}$ and
  $\widehat R:=(\widehat\Delta+1)^{-1}$. Then $\widehat J$ is bounded,
  because
  \begin{equation*}
    \norm{\widehat J}
    =\norm{\wt J J}\leq(1+\delta)(1+\wt\delta)
    \leq 1+ \frac 32(\delta+\wt\delta).
  \end{equation*}
  The second inequality in~\eqref{eq:quasi-uni.a} follows from
  \begin{equation*}
    \norm{\widehat J^*-\widehat J'}
    \leq \norm{J^*(\wt J^*-\wt J')} + \norm{(J^*-J')\wt J'}
    \leq (1+\delta)\wt \delta + \delta(1+2\wt \delta)
    \leq \frac52\delta + \frac52\wt \delta
  \end{equation*}
  as $\norm{\wt J'} \le 1 + 2\delta$ by~\eqref{eq:quasi-uni.op.a''}.
  The first inequality in~\eqref{eq:quasi-uni.b} is also satisfied
  because
  \begin{align*}
    \norm{f - \widehat J' \widehat Jf}
    &\le \norm{f - J'Jf} + \norm{J'(J-J^1)f}
      + \norm{J'(\id_{\wt \HS}-\wt J'\wt J) J^1}
      + \norm{J'\wt J' \wt J(J^1-J)f}\\
    &\le \Bigl(\delta + (1+2\delta)
      \bigl(\delta + \wt \delta (1+\delta)
      + (1+2\wt\delta)(1+\wt\delta)\delta\bigr) \Bigr)\norm[\energy]f
     \le 14(\delta + \wt \delta)\norm[\energy] f
  \end{align*}
  and the second one follows by similar arguments.  Next we prove that
  the two inequalities in~\eqref{eq:quasi-uni.c} also hold. We
  estimate
  \begin{align*}
    \norm{(\widehat J^1-\widehat J)f}
    &\le \norm{(\wt J^1 - \wt J)J^1f} + \norm{\wt J (J^1- J)f}\\
    &\le \bigl(\wt \delta(1+\delta)+(1+\wt \delta)\delta\bigr)\norm[\energy] f
      \le 2(\delta+\wt \delta)\norm[\energy] f
  \end{align*}
  and
  \begin{align*}
    \norm{(\widehat J^{\prime 1}-\widehat J')w}
    &\le \norm{(J^{\prime 1} - J')\wt J^{\prime 1}w}
       + \norm{J'(\wt J^{\prime 1} - \wt J')w}\\
    &\le \bigl(\delta(1+\wt \delta)+(1+2\delta)\wt \delta\bigr)
         \norm[\widehat \energy] w
     \le \frac52 (\delta+\wt \delta)\norm[\widehat \energy] w.
  \end{align*}
  For inequality~\eqref{eq:quasi-uni.d} we estimate
  \begin{align*}
    \bigabs{\widehat \energy(\widehat J^1f,w)-\energy(f,\widehat J^{\prime 1}w)}
    \hspace{-0.2\textwidth}&\\
    &\le \bigabs{\widehat \energy(\wt J^1 J^1f, w)
          - \wt \energy(J^1 f,\wt J^{\prime 1}w)}
    +\bigabs{\wt \energy(\wt J^1f,\wt J^{\prime 1}w)
      - \energy(f, J^{\prime 1} \wt J^{\prime 1}w)}\\
    &\le \wt \delta \norm[\wt \energy] {J^1f} \norm[\widehat \energy] w
      + \delta \norm[\energy] f \norm[\wt \energy] {\wt J^{\prime 1}w}\\
    &\le \bigl(\wt \delta(1+\delta)+\delta(1+\wt \delta)\bigr)
         \norm[\energy] f \norm[\widehat \energy] w
      \le 2(\delta+\wt \delta)\norm[\energy] f \norm[\widehat \energy] w.
      \qedhere
  \end{align*}
\end{proof}

It is a useful feature of \Def{quasi-uni} that it provides us with
some flexibility in terms of the inequalities. The next lemma is one
example.  In~\cite{post-simmer:18} it was applied to avoid a
Poincar\'e-type estimate, i.e., to bypass an estimate of the first
non-zero eigenvalue.
\begin{lemma}[{\cite[Lem.~2.4]{post-simmer:18}}]
  \label{lem:quasi-uni.b}
  Assume that~\eqref{eq:quasi-uni.a} is fulfilled with
  $\deltaA>0$ and~\eqref{eq:quasi-uni.c} with
  $\deltaC>0$. If
  \begin{equation}
    \label{eq:quasi-uni.b'}
    \tag{\ref{eq:quasi-uni.b}'}
    \norm{u - JJ^{\prime1}u}
    \le \delta' \norm[\wt \energy] u \qquad (u \in \wt \HS^1)
  \end{equation}
  holds, then the second inequality in~\eqref{eq:quasi-uni.b} is
  fulfilled with $\delta=\delta'+(1+\deltaA) \deltaC$.

  In particular, if all conditions~\eqref{eq:quasi-uni} are fulfilled
  for some $\delta>0$, except for the second one
  in~\eqref{eq:quasi-uni.b} which is replaced
  by~\eqref{eq:quasi-uni.b'}, then $\energy$ and $\wt \energy$ are
  $\wt \delta$-quasi-unitarily equivalent with $\wt \delta=\delta' +
  (1+\delta)\delta$.
\end{lemma}

\subsection{Consequences of quasi-unitary equivalence}
\label{ssec:que-consequences}


Let $\Delta$ be non-negative and self-adjoint and
$R(z):=(\Delta-z)^{-1}$ be its resolvent.  Let $U$ be an open
neighbourhood of $\spec\Delta \subset \C$ such that $\bd U$ is locally
the graph of a Lipschitz continuous function and such that
$\bd U \cap \spec \Delta = \emptyset$.  Moreover, let $\map \eta U \C$
be a holomorphic function.  Then the integral
\begin{equation}
  \label{eq:def.op.hol}
  \eta(\Delta) := - \frac 1 {2\pi \im} \int_{\bd U} \eta(z)R(z) \dd z
\end{equation}
is defined in the operator norm topology provided
\begin{equation*}
  C_{\eta,\sigma} :=
  \frac 1{2\pi}
  \int_{\bd U} \frac{\abs{\eta(z)}}{d(z,\sigma)} \dd \abs z < \infty
\end{equation*}
for $\sigma:=\spec \Delta$.  For example, if $U$ encloses a compact
subset $K$ of $\spec \Delta$, then $\1_U(\Delta)$ (defined with
$\eta=\1_U$ in~\eqref{eq:def.op.hol}) is the spectral projection onto
$K$.
\begin{theorem}
  \label{thm:func.calc}
  Assume that the forms $\energy$ and $\wt \energy$ corresponding to
  the operators $\Delta$ and $\wt \Delta$ are $\delta$-quasi-unitarily
  equivalent (or that~\eqref{eq:res.z} holds), and that $U$ is an open
  subset such that $\bd U$ is locally Lipschitz and such that
  $\bd U \cap (\spec \Delta \cup \spec{\wt \Delta})=\emptyset$ then
  \begin{equation}
    \label{eq:res.hol}
    \norm{\eta(\wt \Delta)J-J\eta(\Delta)}
    \le C_\eta\delta,
  \end{equation}
  where $C_\eta$ is defined in~\eqref{eq:def.c.eta}.
\end{theorem}
\begin{proof}
  Since the integrals for $\eta(\Delta)$ and $\eta(\wt \Delta)$ exist
  in operator norm, we have
  \begin{equation*}
    \eta(\wt \Delta)J-J\eta(\Delta)
    = -\frac 1{2\pi \im}\int_{\bd U}
    \eta(z)\bigl(\wt R(z)J-J R(z)\bigr) \dd z.
  \end{equation*}
  Taking the operator norm on both sides and using~\eqref{eq:res.z},
  we obtain
  \begin{equation}
    \label{eq:def.c.eta}
    \norm{\eta(\wt \Delta)J-J\eta(\Delta)}
    \le  \underbrace{\frac 1 {2\pi}\int_{\bd U}
        \abs{\eta(z)}
        4 \Bigl(1+\frac{\abs{z+1}}{d(z,\spec \Delta \cup \spec{\wt \Delta})}
        \Bigr)^2 \dd \abs z}_{=:C_\eta} \cdot \delta. \qedhere
  \end{equation}
\end{proof}
Note that $C_\eta<\infty$ implies that $C_{\eta,\spec \Delta}<\infty$
and $C_{\eta,\spec {\wt \Delta}}<\infty$.

\begin{proposition}[{\cite[Lem.~4.2.13]{post:12}}]
  \label{prp:conseq.q-u-e}
  Assume that~\eqref{eq:res.hol} holds.  Then
  \begin{equation*}
    \norm{\eta(\wt \Delta)-J\eta(\Delta)J'}
    \le C'_\eta\delta
    \quadtext{and}
    \norm{\eta(\Delta)-J'\eta(\wt \Delta)J}
    \le C'_\eta \delta
  \end{equation*}
  with
  \begin{equation*}
    C'_\eta:=
    5\sup_{\lambda \in [0,\infty) \cap U} \abs{\eta(\lambda)(\lambda+1)^{1/2}}
      + 3C_\eta
  \end{equation*}
  for all energy forms $\energy$ and $\wt \energy$ (with corresponding
  operators $\Delta$ and $\wt \Delta$, respectively) being
  $\delta$-quasi-unitarily equivalent with identification operators
  $J$ and $J'$ and $\delta \in [0,1]$.
\end{proposition}

Let us calculate explicitly the constants $C_\eta$ and $C'_\eta$ for
two examples of the function $\eta$:
\begin{example}
  \indent
  \begin{enumerate}
  \item \myparagraph{Spectral projections:} Let $I := (a,b)$ such that
    $-1<a<b$ and $a,b \notin \spec \Delta \cup \spec {\wt \Delta}=:S$
    with $d(\{a,b\},S) \ge \eps$ for some $\eps>0$.  We want to
    compare the spectral projections $\1_I(\Delta)$ and
    $\1_I(\wt \Delta)$, defined via the functional calculus for
    self-adjoint operators.  Let
    $U := I \times \im(-\eps,\eps) \subset \C$ be a rectangle
    enclosing $I$.  Note that we have $\1_I(\Delta)=\eta(\Delta)$ with
    $\eta=\1_U$ where the latter operator function is defined via the
    holomorphic functional calculus~\eqref{eq:def.op.hol}; a similar
    statement holds for $\wt \Delta$.  A straightforward estimate
    shows that
    \begin{equation*}
      C_\eta
      = \frac 4\pi \bigl(b-a + \eps\bigr)
      \Bigl(1+\sqrt{1+\Bigl(\frac{b+1}\eps\Bigr)^2}\Bigr)^2
      =\Err(b).
    \end{equation*}
    Moreover, $C_\eta'=5\sqrt{b+1}+3C_\eta=\Err(b)$.

  \item \myparagraph{Heat operator:} For the heat operator, we have
    $\eta_t(\lambda)=\e^{-t\lambda}$ for $t \ge 0$.  Let $U$ be the
    open sector with half-angle $\theta \in (0,\pi/2)$ and vertex at
    $-1$, and symmetric with respect to the real axis.  Then we have
    \begin{equation*}
      C_{\eta_t}
      \le \frac 4 \pi \int_0^\infty
      \e^{-t r \cos \theta}
      \Bigl(1+\frac 1 {\sin\theta}\Bigr)^2 \dd r
      = \frac 4 {\pi\cos \theta} \Bigl(1+\frac 1 {\sin\theta}\Bigr)^2
      \cdot \frac 1t,
    \end{equation*}
    since
    $d(z,\spec \Delta \cup \spec{\wt \Delta})\ge \abs{z+1}\sin\theta$.
    Moreover, as
    \begin{equation*}
      \sup_{\lambda \in [0,\infty)}\abs{\e^{-t\lambda}(\lambda+1)^{1/2}}
      \le
      \begin{cases}
        1, & t \ge 1/2,\\
        1/(2t)^{1/2}, & t \in [0,1/2],
      \end{cases}
    \end{equation*}
    we conclude that
    \begin{equation}
      \label{eq:const.heat.eq}
      C_{\eta_t}'
      =\frac {12} {\pi\cos \theta} \Bigl(1+\frac 1 {\sin\theta}\Bigr)^2
      \cdot \frac 1t + 5.
    \end{equation}
  \end{enumerate}
\end{example}

In particular, we conclude the following convergence result for the
solution of the heat equation:
\begin{corollary}
  \label{cor:conv.sol}
  Let $\energy$ and $\wt {\energy}$ be two $\delta$-quasi-unitarily
  equivalent energy forms with associated operators $\Delta$ and
  $\wt \Delta$.  Assume that $\energy$ and $\wt \energy$ are
  $\delta$-quasi-unitarily equivalent.
   Let $f_t$ resp.\
   $u_t$ be the solution of the heat equations
   \begin{equation*}
     \partial_t f_t + \Delta f_t=0 \quadtext{resp.}
     \partial_t u_t + \wt \Delta u_t=0
   \end{equation*}
   for $t > 0$. If $f_0=J'u_0$, then for any $T>0$ we have
   \begin{equation*}
     \norm[\wt \HS]{u_t - Jf_t} \le C_{\eta_T} \norm[\wt \HS] {u_0}
   \end{equation*}
   for all $t \in [T,\infty)$ with $C'_{\eta_T}=\Err(1/T)$ ($T \to 0$)
   given in~\eqref{eq:const.heat.eq}.
\end{corollary}
\begin{proof}
  We have $u_t=\e^{-t \wt \Delta} u_0$ and
  $f_t = \e^{-t\Delta} J'u_0$.  Then
  \begin{equation*}
    \norm[\wt \HS]{u_t-Jf_t}
    =\norm[\wt \HS]{(\e^{-t \wt \Delta}-J\e^{-t\Delta}J')u_0}
    \le C_{\eta_t}' \delta
  \end{equation*}
  We apply \Prp{conseq.q-u-e} and the concrete estimate for
  $C_{\eta_t}'$ to conclude the desired estimate.
\end{proof}

As in the case of usual norm convergence the operator norm convergence
of spectral projections implies the \emph{convergence of spectra}
(also called \emph{spectral exactness}):
\begin{corollary}[{\cite[Thm.~4.3.3]{post:12}}]
  \label{cor:conv.spec}
  If $\Delta_m$ converges in generalised norm resolvent convergence to
  $\Delta_\infty$, then
  \begin{equation*}
    \spec {\Delta_m} \cap K \to \spec {\Delta_\infty} \cap K
  \end{equation*}
  in Hausdorff distance for any compact $K \subset \R$.
\end{corollary}

If the operators have purely discrete spectrum, we can specify the
error estimate:
\begin{corollary}
  \label{cor:conv.disc.spec}
  Let $\lambda_k(\Delta_m)$ resp. $\lambda_k(\Delta_\infty)$ denote
  the $k^{\text{th}}$ eigenvalue of $\Delta_m$ resp. $\Delta_\infty$
  (in increasing order and with respect to multiplicity). Then
  \begin{equation*}
    \bigabs{\lambda_k(\Delta_m)-\lambda_k(\Delta_\infty)}\leq C_k\delta_m
  \end{equation*}
  for all $m\in\N$ such that $\dim\HS_m\geq k$, where $C_k$ depends
  only on $\lambda_k(\Delta_\infty)$.
\end{corollary}

In the case of purely discrete spectrum (or isolated eigenvalues) we
can approximate an eigenfunction also in energy norm:

\begin{proposition}[{\cite[Prp.~2.6]{post-simmer:18}}]
  \label{prp:conv.ef}
  Let $\energy$ and $\wt {\energy}$ be two $\delta$-quasi-unitarily
  equivalent energy forms with associated operators $\Delta$ and $\wt
  \Delta$.  Assume that $\wt \Phi$ is an eigenvector of $\wt \Delta$,
  such that its eigenvalue $\wt \lambda$ is discrete in $\spec {\wt
    \Delta}$, i.e., there is an open disc $D$ in $\C$ such that $\spec
  {\wt \Delta} \cap D=\{\wt \lambda\}$.  Then there exists a
  normalised eigenvector $\Phi$ of $\Delta$ with $\Phi \in \ran
  \1_D(\Delta)$ and a universal constant $C$ depending only on $\wt
  \lambda$ (and the radius of $D$) such that
  \begin{equation*}
    \norm[\wt \energy]{J^1 \Phi - \wt \Phi} \le C\delta.
  \end{equation*}
\end{proposition}

Note that the eigenvalue $\wt \lambda$ does not necessarily need to
have finite multiplicity.

%
\section{Post-critically finite self-similar fractals}
\label{sec:pcf-fractals}
%

In~\cite{post-simmer:18} the authors applied the quasi-unitary
equivalence to the case of certain fractals called
\emph{post-critically finite self-similar sets} (which supports a
regular resistance form in the sense
of~\cite{kigami:01}; see also~\cite{strichartz:06}). Here, we will simply discuss two
examples. For the general case, we refer to~\cite{post-simmer:18} and
for a further generalisation to magnetic energy forms on finitely
ramified fractals, we refer to~\cite{post-simmer:pre18c}.


\subsection{The unit interval}
\label{ssec:ex.unit-interval}


At the first glance it might look a bit odd to call the unit interval
$K=[0,1]$ a fractal, but it will turn out that this approach is quite
elegant for the approximation. We begin by defining two contractions
$F_1,F_2\colon\R\to\R$ with contraction ratio $1/2$ and fixed points
$0$ and $1$ by
\begin{equation*}
F_1(t)=\frac t 2
\qquad\text{and}\qquad
F_2(t)=\frac {1+t}2.
\end{equation*}
Then, we have $K=F_1(K)\cap F_2(K)$ and $K$ is the unique non-empty
compact subset of $\R$ with that property. We call $K$ the
\emph{self-similar fractal} with respect to $F=\{F_1,F_2\}$. Moreover,
the maps $F_j$ describe a cell structure on $K$ via
\begin{equation*}
  w\mapsto F_w(K):=(F_{w_1}\circ\dots\circ F_{w_m})(K)
\end{equation*}
where $w=w_1\dots w_m\in W_m:=\{1,2\}^m$ is a word of length $m$. We
refer to $F_w(K)$ as an $m$-cell whenever $w\in W_m$.

Next, we define the \emph{(vertex) boundary} by $V_0=\{0,1\}$. Note
that in the special case of the interval, the topological boundary and
the vertex boundary coincide but that is not necessarily the case (see
e.g.~the Sierpi\'nski gasket). Then, we define the approximating
sequence of (finite weighted discrete) graphs as follows: Let
$G_0=(V_0,E_0)$ be the complete graph with two vertices $V_0=\{0,1\}$
and one edge: Moreover, define $G_m=(V_m,E_m)$ inductively, where
$V_m=\set{k2^{-m}}{k=0,\dots,2^m}$ are the $m$-dyadic numbers and
where we have an edge between (distinct) vertices $x$ and $y$ in $V_m$
if and only if $\abs{x-y}=2^{-m}$.

Let us specify the Hilbert spaces and energy forms now. As a measure
on $K=[0,1]$ we fix the Lebesgue measure $\mu$ and our Hilbert space
is the usual space of square integrable function with respect to the
Lebesgue measure, i.e., $\wt\HS=\Lsqr{K,\mu}$.

The approximating measure $\mu_m=\{\mu_m(x)\}_{x\in V_m}$ on $G_m$ is defined by
\begin{equation}
  \label{eq:approx-measure-interval}
  \mu_m(x)=\int_0^1\psi_{x,m}(t)\dd\mu(t)=\frac 1 {2^m}
  \begin{cases}
    1\quad & x\in V_m\setminus V_0\\
    1/2\quad & x\in V_0,
  \end{cases}
\end{equation}
where $\psi_{x,m}\colon K\to [0,1]$ is given by $\1_{\{x\}}$ on $V_m$
and extended to $K$ by linear interpolation. Hence, our Hilbert space
on the graph $G_m$ is $\HS_m=\lsqr{V_m,\mu_m}$ with norm
\begin{equation*}
  \normsqr[\lsqr{V_m,\mu_m}]{f}=\sum_{x\in V_m}\mu_m(x)\abssqr{f(x)}.
\end{equation*}

On each graph $G_m$ we now define a discrete energy form $\energy_m$
in $\lsqr{V_m,\mu_m}$. For each
$f\in\ell(V_m):=\set{f}{f\colon V_m\to\C}$, we set
\begin{equation}
  \label{eq:disc-energy-interval}
  \energy_m(f)=\sum_{\{x,y\} \in E_m}c_{\{x,y\},m}\abssqr{f(y)-f(x)},
\end{equation}
where the \emph{conductances} $c_{\{x,y\},m}\geq 0$ are chosen such
that
\begin{equation}
  \label{eq:comp}
  \energy_m(\phi)
  =\min\bigset{\energy_{m+1}(f)}{f\in\ell(V_{m+1}),f\restr{V_m}=\phi},
\end{equation}
for all $\phi\in\ell(V_m)$. Working this out, we see that
$c_{\{x,y\},m}=2^m$. A sequence $\{\energy_m\}_{m\in\N_0}$ of energy
forms that satisfies~\eqref{eq:comp} for all $m$ is called
\emph{compatible sequence}.

From the classical theory of calculus it is well-known that the limit
form is given by
\begin{equation*}
  \energy(u)=\int_0^1\abssqr{u^\prime(t)}\dd\mu(t)
\end{equation*}
for each weakly differentiable $u \in\Lsqr{K,\mu}$ with
$u' \in\Lsqr{K,\mu}$, i.e., we have $\dom\energy=\Sob{K,\mu}$.

\begin{theorem}[{\cite{post-simmer:18}}]
  The energy form $(\energy,\Sob{K,\mu})$ in $\Lsqr{K,\mu}$ and the
  discrete energy form $\energy_m$ in $\lsqr{V_m,\mu_m}$ are
  $\delta_m$-quasi-unitarily equivalent, where the error is
  \begin{equation*}
    \delta_m=(1+\sqrt{2})\cdot\frac 1 {2^m}.
  \end{equation*}
\end{theorem}


Let us briefly discuss the idea of the proof: First, we need to choose
the identification operators from \Def{quasi-uni}. On the Hilbert
space level, we define $J_m\colon\HS_m\to\wt\HS$ and
$J_m^\prime=J_m^\ast\colon\wt\HS\to\HS_m$ by
\begin{equation*}
  J_mf:=\sum_{x\in V_m}f(x)\psi_{x,m}
  \qquad
  \text{resp.}
  \qquad
  J_m^\prime u(y)=\frac 1 {\mu_m(y)}\iprod[\wt\HS]{u}{\psi_{y,m}},
\end{equation*}
where $f\in\HS_m$, $u\in\wt\HS$ and $y\in V_m$. Moreover, we define
$J_m^1\colon\HS_m^1\to\wt\HS^1$ by $J_m^1f:=J_mf$. Note that this is
well-defined because $\psi_{x,m}\in\dom\energy$. The last operator
$J_m^{\prime 1}$ is chosen to be the evaluation in points of $V_m$,
i.e., $J_m^{\prime 1}u(x)=u(x)$. Again, this choice makes sense
because functions in the domain of $\energy$ are continuous on
$K$. That is because
\begin{equation}
  \label{eq:hoelder-estimate}
  \abssqr{u(x)-u(y)}
  \leq \energy(u)R(x,y),
\end{equation}
for $x,y\in K$, where $R$ is the \emph{resistance metric} associated
with $\energy$, given by
\begin{equation*}
  R(x,y)
  :=\Bigl(\min\bigset{\energy(u)}{u\in\dom\energy, u(x)=1
    \text{ and } u(y)=0}\Bigl)^{-1}
\end{equation*}
and since $R(x,y)=\abs{x-y}$, the relative topology of $K$ coincides
with the $R$-topology (see~\cite[Sec.~3.4]{kigami:01} for a more
general result).

Now, we need to verify the validity of the inequalities in
\Def{quasi-uni}. This is done in~\cite[Sec.~4]{post-simmer:18} in
greater details but let us discuss the key steps here:

Applying the Cauchy-Young inequality in the first inequality
of~\eqref{eq:quasi-uni.a} we see that
$\norm[\wt\HS]{Jf}\leq\norm[\HS_m]{f}$ for each $f\in\HS_m$ and the
second one is fulfilled because $J_m^\prime=J_m^\ast$.

The inequalities in~\eqref{eq:quasi-uni.b} follow by applying the
Cauchy-Schwarz inequality and by using the improved H\"older
inequality~\eqref{eq:hoelder-estimate}. For the first one, we rewrite
\begin{equation*}
f(y)=\frac 1 {\mu_m(y)}\sum_{x\in V_m}f(x)\iprod[\wt\HS]{\psi_{x,m}}{\psi_{y,m}}
\end{equation*}
using the fact that $\{\psi_{x,m}\}_{x\in V_m}$ is a partition of unity on $K$ and
\begin{equation*}
J_m^\prime J_mf(y)
=\sum_{x\in V_m} f(x) J_m^\prime \psi_{x,m}(y)
=\frac 1 {\mu_m(y)}\sum_{x\in V_m} f(x)\iprod[\wt\HS]{\psi_{x,m}}{\psi_{y,m}}.
\end{equation*}
Hence, by applying the above mentioned inequalities and some standard
arguments, we can estimate $f-J_m^\prime J_mf$ in norm.

Note, that the first inequality from~\eqref{eq:quasi-uni.c} is
trivially fulfilled by the choice of $J_m$ and $J_m^1$. Instead of
verifying the second one, we apply \Lem{quasi-uni.b}. This is
particularly useful here because it helps us to skip a discussion
about eigenvalues, we would otherwise have (see
e.g.~\cite{post-simmer:21}).

The particular choice of the identification operators becomes clear
now: The last inequality~\eqref{eq:quasi-uni.d} holds actually with
equality because the $\{\psi_{x,m}\}$ minimise the energy, i.e. as
above,
$\energy(\psi_{x,m})=\energy_m(\psi_{x,m}\restr{V_m})=\energy_m(\1_{\{x\}})$
where $\1_{\{x\}}$ is the characteristic function of the set
$\{x\}\subset V_m$. Note that the letter expression can be computed
explicitly using~\eqref{eq:disc-energy-interval}.


\subsection{The Sierpi\'nski gasket}
\label{ssec:ex.SG}


A more illustrative example for a post-critically finite self-similar
fractal is the Sierpi\'nski gasket which is described by the family of
contractions $F$, given by
\begin{align*}
  F_j\colon\R^2\to\R^2,
  && F_j(x)=\frac 1 2 (x-p_j)+p_j
  && (j=1,2,3)
\end{align*}
where the fixed points $p_j$ are chosen such that $\{p_1,p_2,p_3\}$
are the vertices of an equilateral triangle in $\R^2$. Then, as in the
case of the unit interval, the Sierpi\'nski gasket is defined as the
unique non-empty compact subset $K$ of $\R^2$ that satisfies
\begin{equation*}
K=F(K):=F_1(K)\cup F_2(K)\cup F_3(K).
\end{equation*}
Again, the family of contractions $\{F_1,F_2,F_3\}$ describes a cell
structure on the Sierpi\'nski gasket via the map $w\mapsto F_w(K)$
where $w\in W_m=\{1,2,3\}^m$. The vertex boundary is defined as
$V_0=\{p_1,p_2,p_3\}$. Note that in contrast to the situation in the
example of the unit interval, $V_0$ does not coincide with the
topological boundary of $K$ which is actually $K$ itself.

We define our approximating sequence of graphs in the same way as
before: Let $G_0=(V_0,E_0)$ be the complete graph and let
$G_m=(V_m,E_m)$ be given by
\begin{align*}
V_m:=\bigcup_{j=1}^3F_j(V_{m-1})
&& \text{and} &&
E_m:=\bigset{e}{e=\text{ and }x\sim_my},
\end{align*}
where we write $x\sim_my$ if and only if $x$ and $y$ are two distinct
vertices in $V_m$ and there exists a word $w\in W_m$ such that
$x,y\in F_w(K)$. Moreover, we define an energy form on $G_m$ by
\begin{equation*}
\energy_m(f)=\sum_{x\sim_my}\Bigl(\frac 5 3 \Bigr)^m\abssqr{f(x)-f(y)},
\end{equation*}
for each $f\in\ell(V_m)$. Here we sum over all vertices $x\in V_m$ and
their neighbours $y$ in the level $m$ graph and, as in the case of the
interval, the conductances $c_{\{x,y\},m}=(5/3)^m$ are chosen such
that the sequence of energy forms $\{\energy_m\}_{m\in\N_0}$ is
compatible.

Then the limit form exists and we define an energy form on the
Sierpi\'nski gasket by
\begin{align*}
\energy(u):=\lim_{m\to\infty}\energy_m(u\restr{V_m}),
&&
u\in\dom\energy:=\bigset{u}{u\colon K\to\C,\energy(u)<\infty}
\end{align*}
(see~\cite{kigami:01,strichartz:06}).


As a (canonical) measure $\mu$ on the Sierpi\'nski gasket we choose
the homogeneous self-similar measure, i.e., the uniquely determined
probability measure $\mu$ that satisfies
\begin{equation*}
\mu=\frac 1 3 \bigl(\mu\circ F_1^{-1}+\mu\circ F_2^{-1}+\mu\circ F_3^{-1}\bigr).
\end{equation*}
That is, $\mu$ is the Hausdorff measure of dimension $\log 3/\log 2$
and every $m$-cell $F_w(K)$ has measure $1/3^m$ (we would like to
stress that our approach works for a general Borel regular probability
measure on $K$; see~\cite{post-simmer:18} for details). As Hilbert
space structure on the fractal $K$, we choose
$\wt\HS=\Lsqr{K,\mu}$. Then $(\energy,\dom\energy)$ is a closed
quadratic form in $\wt\HS$.

On the graphs $G_m$ we define a measure as
in~\eqref{eq:approx-measure-interval} but here we choose the functions
$\psi_{x,m}\colon K\to[0,1]$ to be the unique solution of the
minimisation problem
\begin{equation*}
  \energy_m(\1_{\{x\}})
  =\min\bigset{\energy(u)}{u\in\dom\energy,u\restr{V_m}=\1_{\{x\}}}.
\end{equation*}
These functions exist and are called $m$-harmonic functions with
boundary values $\1_{\{x\}}$ on $V_m$ (cf.~\cite{kigami:01}). The
values of $\psi_{x,m}$ can be computed explicitly by iteration: If the
values in the vertices of $V_m$ are known, then, for each vertex
$y\in V_{m+1}\setminus V_m$, there exists a unique $m$-cell that
contains $y$; the value $\psi_{x,m}(y)$ is given by $1/5$ times the
value at the vertex in $V_m$ opposite to $y$ in the $m$-cell
plus $2/5$ times the values of $\psi_{x,m}$ at the vertices (of $V_m$)
adjacent to $y$ in the same $m$-cell
(cf.~\cite[Sec.~1.3]{strichartz:06}).

By the symmetry of the Sierpi\'nski gasket and the functions
$\{\psi_{x,m}\}_{x\in V_m}$, which define a partition of unity on $K$,
we can specify $\mu_m$ also in this example as
\begin{equation*}
\mu_m(x)=\int_K\psi_{x,m}\dd\mu=\frac 1 {3^m}
\begin{cases}
1/3 \qquad &x\in V_0\\
2/3 &x\in V_m\setminus V_0.
\end{cases}
\end{equation*}
The Hilbert space, we consider on the approximating sequence of graphs
is again given by $\HS_m=\lsqr{V_m,\mu_m}$ and we conclude:


\begin{theorem}[{\cite{post-simmer:18}}]
  The energy form $(\energy,\dom\energy)$ in $\Lsqr{K,\mu}$ and the
  discrete energy form $\energy_m$ in $\lsqr{V_m,\mu_m}$ are
  $\delta_m$-quasi-unitarily equivalent where the error is
  \begin{equation*}
    \delta_m=\frac{(1+\sqrt{3})\sqrt{2}}{\sqrt 3 }\cdot\frac 1 {5^{m/2}}.
  \end{equation*}
\end{theorem}


The idea of the proof is the same as described above in the case of
the unit interval (see~\cite{post-simmer:18}).


%
\section{Neumann obstacles}
\label{sec:neu-obstacles}
%
In this section we briefly present another class of examples.  For
details, we refer to~\cite{anne-post:21}.  Let $X$ be a complete
Riemannian manifold of bounded geometry (i.e., its Ricci curvature is
bounded from below and its injectivity radius is bounded from below by
a positive constant); a simple example is $X=\R^n$.  We denote by
$\Lsqr X$ the Hilbert space of square-integrable functions with
respect to the standard volume measure, and by $\Sob[k] X$ the Sobolev
space of $k$-times weakly differentiable and square integrable
functions.  Denote its Laplacian by $\Delta_X \ge 0$ (defined via its
quadratic form $\energy_X(f)=\int_X \abssqr{df} \dvol$).  Under the
above assumptions (completeness and bounded geometry) it can then be
shown that there is a constant $\Cellreg \ge 1$ such that
\begin{equation}
  \label{eq:ell.reg}
  \norm[{\Sob[2]X}] f \le \Cellreg \norm[\Lsqr X]{(\Delta_X+1)f}
\end{equation}
for all $f \in \dom \Delta_X=\Sob[2] X$.

We assume that $B \subset X$ is a closed subset such that the
following holds:
\begin{enumerate}
\item there is $\delta \ge 0$ such that
  \begin{equation*}
    \norm[\Lsqr B] f \le \delta \norm[\Sob X] f
  \end{equation*}
  for all $f \in \Sob X$;
\item there is a \emph{bounded extension operator}, i.e., there is
  $\map E {\Sob{X \setminus B}}{\Sob X}$ such that
    $E u \restr {X \setminus B} = u$ with operator norm bounded by
    $\Cext \ge 1$.
\end{enumerate}
One can think of $B$ as the disjoint union of small balls or other
obstacles.  Denote by $\Delta^\Neu_{X \setminus B}$ the Neumann
Laplacian defined via its quadratic form
$\energy^\Neu_{X \setminus B}(u):=\int_{X \setminus B} \abssqr{du}
\dvol$. It can be seen that the first estimate extends to
\begin{equation*}
    \norm[\Lsqr B] {df} \le \delta \norm[{\Sob[2]X}] f
\end{equation*}
for $f \in \Sob[2] X$ without any assumption on the manifold
(cf.~\cite[Prp.~3.7]{anne-post:21}).  We have the following result:
\begin{theorem}[{\cite[Thm.~4.3]{anne-post:21}}] Under the above assumptions,
  the Laplacian $\Delta_X$ and the Neumann Laplacian
  $\Delta^\Neu_{X\setminus B}$ are $\Cext\Cellreg \delta$-quasi-unitarily
  equivalent.
\end{theorem}
\begin{proof}
  We are showing a slightly modified version of quasi-unitary
  equivalence for the corresponding energy forms.  We first define the
  following identification operators as follows:
  \begin{equation*}
    \map J {\HS:=\Lsqr X}{\wt \HS:=\Lsqr{X \setminus B}},
    \qquad
    f \mapsto f \restr {X \setminus B},
  \end{equation*}
  $\map {J^1}{\HS^1:=\Sob X}{\wt \HS^1:=\Sob {X \setminus B}}$,
  $f \mapsto f \restr{X \setminus B}$, $J'=J^*$ (hence $J'u$ is the
  extension of $u$ by $0$ onto $B$) and
  \begin{equation*}
    \map {J^{\prime 1}}{\wt \HS:=\Sob{X \setminus B}}{\HS^1=\Sob X},
    \qquad
    u \mapsto Eu.
  \end{equation*}
  In particular, $J^1f=Jf$ for $f \in \Sob X$, $JJ'u=u$, $\norm J=1$
  and $J'=J^*$.  It remains to check the first inequality
  of~\eqref{eq:quasi-uni.b}, the second of~\eqref{eq:quasi-uni.c} and
  a modified version of~\eqref{eq:quasi-uni.d}: the first estimate is
  fulfilled with $\delta$ since
  \begin{equation*}
    \norm[\Lsqr X]{f-J'Jf}
    =\norm[\Lsqr B] f
    \le \delta \norm[\Sob X] f
  \end{equation*}
  by our assumption.  For the second, we argue
  \begin{equation*}
    \norm[\Lsqr X]{J^{\prime 1}u-J'u}
    =\norm[\Lsqr B] {Eu}
    \le \delta \norm[\Sob X] {Eu}
    \le \Cext \delta \norm[\Sob {X \setminus B}] u,
  \end{equation*}
  hence the estimate is fulfilled with $\Cext \delta \ge \delta$.  Instead
  of~\eqref{eq:quasi-uni.d}, we show the slightly stronger estimate
  \begin{equation}
    \label{eq:quasi-uni.d2}
    \bigabs{\wt \energy(J^1f, u) - \energy(f, J^{\prime1}u)}
    \le \delta \norm[\Lsqr X] {(\Delta_X+1)f} \norm[\wt \energy] u
    (f \in \HS^2:=\dom \Delta_, u \in \wt \HS^1),
  \end{equation}
  i.e., we use the operator graph norm instead of the energy norm for
  $f$.  In particular, we have
  \begin{align*}
    \bigabs{\wt \energy(J^1f, u) - \energy(f, J^{\prime1}u)}
    &= \bigabs{\int_B \iprod {df}{d(Eu)} \dvol}\\
    &\le \norm[\Lsqr B]{df} \norm[\Lsqr B]{d(Eu)}
    \le \norm[{\Sob[2] X}] f \norm[\Sob X]{Eu}\\
    &\le \Cellreg \norm[\Lsqr X]{(\Delta_X+1)f} \Cext\norm[\Sob{X\setminus B}] u
  \end{align*}
  using the assumptions. The quasi-unitary equivalence for the
  operators follows then similarly as in \Prp{q-u-e.qf.op}.
\end{proof}
Note that if $B=\bigcup_{x \in I_\eps} B_\eps(x)$ is the disjoint
union of balls of radius $\eps>0$ with centres $x \in I_\eps$
separated by $2\eps^\alpha$ (i.e., $x,y\in I_\eps$ and $x \ne y$,
implies $d(x,y) > 2 \eps^\alpha$) with $0<\alpha<1$, then one can show
that the extension operator is \emph{uniformly} bounded, i.e., $\Cext$
can be chosen to be independent of $\eps$.  Moreover, one can choose
$\delta=\delta_\eps$ to be of order $\eps^{1-\alpha}$ in dimension
$n \ge 3$ (resp.\ $\eps^{1-\alpha}\log\abs \eps$ in dimension $n=2$).
Note that the sets $I_\eps$ for different values of $\eps$ may be
totally unrelated, see~\cite[Sec.~4.2]{anne-post:21} for details.

%
%


\providecommand{\bysame}{\leavevmode\hbox to3em{\hrulefill}\thinspace}
\providecommand{\MR}{\relax\ifhmode\unskip\space\fi MR }
\providecommand{\MRhref}[2]{%
  \href{http://www.ams.org/mathscinet-getitem?mr=#1}{#2}
}
\providecommand{\href}[2]{#2}


\end{document}

%% file: olaf-ams-layout-macros.tex



\newcommand{\myfont}{\sffamily}

\newtheoremstyle{mythmstyle}
  {\topsep}
  {\topsep}
  {\itshape}
  {}
  {\bfseries \myfont}
  {.}
  {.5em}
  {}

\newtheoremstyle{mydefstyle}
  {\topsep}
  {\topsep}
  {\normalfont}
  {}
  {\bfseries \myfont}
  {.}
  {.5em}
  {}

\swapnumbers                    

\theoremstyle{mythmstyle}       






\let\MakeUppercase\relax

\expandafter\let\expandafter\oldproof\csname\string\proof\endcsname
\let\oldendproof\endproof
\renewenvironment{proof}[1][\bfseries\myfont\proofname]{%
  \oldproof[\bfseries \myfont #1]%
}{\oldendproof}


\makeatletter
\newtoks\thm@headfont  \thm@headfont{\bfseries \myfont}

\renewcommand\section{\@startsection{section}{1}%
  \z@{.7\linespacing\@plus\linespacing}{.5\linespacing}%
  {\Large\myfont\bfseries}}
\renewcommand\subsection{\@startsection{subsection}{2}%
  \z@{-.5\linespacing\@plus-.7\linespacing}{.5\linespacing}%
  {\large\myfont\bfseries}}
\renewcommand\subsubsection{\@startsection{subsubsection}{3}%
  \z@{.5\linespacing\@plus.7\linespacing}{-.5em}%
  {\myfont\bfseries}}
\renewenvironment{abstract}{%
  \ifx\maketitle\relax
    \ClassWarning{\@classname}{Abstract should precede
      \protect\maketitle\space in AMS document classes; reported}%
  \fi
  \global\setbox\abstractbox=\vtop \bgroup
    \normalfont\Small
    \list{}{\labelwidth\z@
      \leftmargin3pc \rightmargin\leftmargin
      \listparindent\normalparindent \itemindent\z@
      \parsep\z@ \@plus\p@
      
    }%
    \item[\hskip\labelsep
      \myfont\bfseries
    \abstractname.]%
}{%
  \endlist\egroup
  \ifx\@setabstract\relax \@setabstracta \fi
}
\renewcommand\contentsnamefont{\myfont\bfseries}
\renewcommand\@starttoc[2]{\begingroup
  \setTrue{#1}%
  \par\removelastskip\vskip\z@skip
  \@startsection{}\@M\z@{\linespacing\@plus\linespacing}%
    {.5\linespacing}{
      \contentsnamefont}{#2}%
  \ifx\contentsname#2%
  \else \addcontentsline{toc}{section}{#2}\fi
  \makeatletter
  \@input{\jobname.#1}%
  \if@filesw
    \@xp\newwrite\csname tf@#1\endcsname
    \immediate\@xp\openout\csname tf@#1\endcsname \jobname.#1\relax
  \fi
  \global\@nobreakfalse \endgroup
  \addvspace{32\p@\@plus14\p@}%
  \let\tableofcontents\relax
}
\renewcommand\@settitle{\begin{center}%
  \baselineskip14\p@\relax
    \LARGE
    \bfseries
    \myfont
  \@title
  \end{center}%
}
\renewcommand\@setauthors{%
  \begingroup
  \def\thanks{\protect\thanks@warning}%
  \trivlist
  \centering\footnotesize \@topsep30\p@\relax
  \advance\@topsep by -\baselineskip
  \item\relax
  \author@andify\authors
  \def\\{\protect\linebreak}%
  \large
  \myfont\bfseries\authors
  \ifx\@empty\contribs
  \else
    ,\penalty-3 \space \@setcontribs
    \@closetoccontribs
  \fi
  \endtrivlist
  \normalfont\myfont\@setaddresses
  \endgroup
}
\renewcommand\@setaddresses{\par
  \nobreak \begingroup
\footnotesize
  \def\author##1{\nobreak\addvspace\bigskipamount}%
  \def\\{\unskip, \ignorespaces}%
  \interlinepenalty\@M
  \def\address##1##2{\begingroup
    \par\addvspace\bigskipamount\indent
    \@ifnotempty{##1}{(\ignorespaces##1\unskip) }%
    {
      \ignorespaces##2}\par\endgroup}%
  \def\curraddr##1##2{\begingroup
    \@ifnotempty{##2}{\nobreak\indent\curraddrname
      \@ifnotempty{##1}{, \ignorespaces##1\unskip}\/:\space
      ##2\par}\endgroup}%
  \def\email##1##2{\begingroup
    \@ifnotempty{##2}{\nobreak\indent\emailaddrname
      \@ifnotempty{##1}{, \ignorespaces##1\unskip}\/:\space
      \ttfamily##2\par}\endgroup}%
  \def\urladdr##1##2{\begingroup
    \def~{\char`\~}%
    \@ifnotempty{##2}{\nobreak\indent\urladdrname
      \@ifnotempty{##1}{, \ignorespaces##1\unskip}\/:\space
      \ttfamily##2\par}\endgroup}%
  \addresses
  \endgroup
}
\renewcommand\enddoc@text{\ifx\@empty\@translators \else\@settranslators\fi
}
\renewcommand\@secnumfont{\myfont\bfseries} 

\renewcommand\maketitle{\par
  \@topnum\z@ 
  \@setcopyright
  \thispagestyle{firstpage}
  \ifx\@empty\shortauthors \let\shortauthors\shorttitle
  \else \andify\shortauthors
  \fi
  \@maketitle@hook
  \begingroup
  \@maketitle
  \toks@\@xp{\shortauthors}\@temptokena\@xp{\shorttitle}%
  \toks4{\def\\{ \ignorespaces}}
  \edef\@tempa{%
    \@nx\markboth{\the\toks4
      \@nx\MakeUppercase{\the\toks@}}{\the\@temptokena}}%
  \@tempa
  \endgroup
  \c@footnote\z@
  \@cleartopmattertags
}
\def\@captionheadfont{\myfont\bfseries} 

\makeatother

%% file: olaf-thm-environment-macros.tex
%
\usepackage{amsthm}


\swapnumbers  
\newcounter{intro}

\swapnumbers 
\newtheorem{theorem}{Theorem}[section]
\newtheorem{proposition}[theorem]{Proposition}
\newtheorem{lemma}[theorem]{Lemma}
\newtheorem{corollary}[theorem]{Corollary}

\theoremstyle{mydefstyle}        
\newtheorem{definition}[theorem]{Definition}

\newtheorem{example}[theorem]{Example}


\newtheorem*{remark*}{Remark}


\numberwithin{equation}{section}

%% file: olaf-refs-unified.tex
%
\newcommand{\Sec}[1]{Section~\ref{sec:#1}}

\newcommand{\Subsec}[1]{Subsection~\ref{ssec:#1}}




\newcommand{\Thm}[1]{Theorem~\ref{thm:#1}}

\newcommand{\Lem}[1]{Lemma~\ref{lem:#1}}

\newcommand{\Prp}[1]{Proposition~\ref{prp:#1}}

\newcommand{\Def}[1]{Definition~\ref{def:#1}}

%% file: olaf-math-symbols.tex
%
%
\usepackage{accents}
\newcommand{\abs}[2][{}]{\lvert{#2}\rvert_{{#1}}}    
\newcommand{\abssqr}[2][{}]{\lvert{#2}\rvert^2_{#1}} 
\newcommand{\bigabs}[2][{}]{\bigl\lvert{#2}\bigr\rvert_{#1}}     

\newcommand{\normsymb}{\|}
\newcommand{\bignormsymb}[1]{#1\|}

\newcommand{\norm}[2][{}]{\normsymb{#2}\normsymb_{{#1}}}    
\newcommand{\normsqr}[2][{}]{\normsymb{#2}\normsymb^2_{#1}} 
\newcommand{\bignorm}[2][{}]{\bignormsymb{\bigl}{#2}\bignormsymb{\bigr}_{#1}}



\newcommand{\iprod}[3][{}]{\langle{#2},{#3}\rangle_{#1}}  
\newcommand{\bigiprod}[3][{}]{\bigl\langle{#2},{#3}\bigr\rangle_{#1}}

\newcommand{\set}[2]{\{ \, #1 \, | \, #2 \, \} }      
\newcommand{\bigset}[2]{\bigl\{ \, #1 \, \bigl|\bigr. \, #2 \, \bigr\} }

\newcommand{\map}[3]{ #1 \colon #2 \longrightarrow #3}    

\newcommand{\bd}  {\partial}          

\newcommand{\restr}[1]{{\restriction}_{#1}} 

%

\def\XXint#1#2#3{{\setbox0=\hbox{$#1{#2#3}{\int}$}
     \vcenter{\hbox{$#2#3$}}\kern-.5\wd0}}
\def\XXsum#1#2#3{{\setbox0=\hbox{$#1{#2#3}{\int}$}
     \vcenter{\hbox{$#2#3$}}\kern-.60\wd0}}




\DeclareMathOperator{\dd}    {d\!}

\DeclareMathOperator{\dom}    {dom}
\DeclareMathOperator{\ran}    {ran}
\DeclareMathOperator{\id}     {id}   

\DeclareMathOperator{\dvol}    {d\, vol}

\newcommand{\specsymb} {\sigma} 

\newcommand{\spec}[2][{}]   {\specsymb_{\mathrm{#1}}(#2)}


\newcommand{\eps}{\varepsilon} 

\renewcommand{\phi}{\varphi}   
\renewcommand{\rho}{\varrho}   

\DeclareMathOperator{\myRe} {Re}
\renewcommand{\Re}     {\myRe}
\DeclareMathOperator{\myIm}     {Im}
\renewcommand{\Im}     {\myIm}

\newcommand{\conj}[1]{\overline {#1}}

\newcommand{\R}{\mathbb{R}} 
\newcommand{\C}{\mathbb{C}} 
\newcommand{\N}{\mathbb{N}} 

\newcommand{\1}{\mathbbm 1}                    

\newcommand{\e}{\mathrm e}  
\newcommand{\im}{\mathrm i} 


\newcommand{\wt}{\widetilde}           

\newcommand{\energy}{\mathcal E}  


\newcommand{\Neu}{{\mathrm N}}              

\newcommand{\Err}{\mathrm O}



%% file: olaf-other-macros-myenumi-myparagraph.tex
\newcounter{myenumi}



\newcommand{\myparagraph}[1]{\noindent\textbf{\myfont{#1}}}

\newcommand{\quadtext}[1]{\quad\text{#1}\quad}

%% file: olaf-spaces-macros.tex
\newcommand{\HS}{\mathscr H}           


\newcommand{\Sobsymb} {\mathsf H} 

\newcommand{\Lsymb}    {\mathsf L}     
\newcommand{\lsymb}    {\ell}          

\newcommand{\Sobspace}[1][1]{\Sobsymb^{{#1}}}

\newcommand{\Lpspace}[1][p]    {\Lsymb_{#1}}     
\newcommand{\lpspace}[1][p]    {\lsymb_{#1}}     
\newcommand{\Lsqrspace}    {\Lpspace[2]}     
\newcommand{\lsqrspace}    {\lpspace[2]}          





\newcommand{\Lsqr}[2][{}]{\Lsqrspace^{#1}({#2})} 
 %

\newcommand{\lsqr}[2][{}]{\lsqrspace^{#1}({#2})}   



\newcommand{\Sob}[2][1]{\Sobspace [{#1}]({#2})}         















